\documentclass[12pt,a4paper]{article}
\RequirePackage{etex}
\usepackage{amsfonts, amsmath, amssymb,amsthm,graphics,mathrsfs,verbatim, CJK}
\usepackage{makeidx}
\usepackage{palatino}
\usepackage{titletoc}
\usepackage{tocloft}
\usepackage{titlesec}
\usepackage{stmaryrd,color}
\usepackage[colorlinks=true]{hyperref}
\usepackage{graphicx}
\usepackage{pgfplots}
\usepackage{tikz}
\usepackage{pstricks,pst-node,multido,ifthen,calc}

\titleformat{\section}{\Large}{}{0.5em}{}
\makeindex

\normalsize

\def\nin{\noindent}

\def\seq{\subseteq}

\def\fs{\footnotesize}

\def\vs{\vspace*}

\newtheorem{mthm}{Theorem}[section]
\newtheorem{mylem}[mthm]{Lemma}
\newtheorem{myprn}[mthm]{Proposition}
\newtheorem{mycor}[mthm]{Corollary}
\newtheorem{mydef}[mthm]{Definition}
\newtheorem{myrem}[mthm]{Remark}
\newtheorem{mycon}[mthm]{Construction}
\newtheorem{myeg} [mthm]{Example}
\newtheorem{myque} [mthm]{Question}
\newtheorem{myalg} [mthm]{Algorithm}

\newenvironment{prof}{\noindent $Proof.$ \rm}{\hfill $\Box$}

\def\fs{\footnotesize}

\def\seq{\subseteq}

\def\vs{\vspace*}

\def \nin {\noindent}

\def \Lemma #1 {\vs{2mm}\nin {\bf Lemma #1} }
\def \Prop #1 {\vs{2mm}\nin {\bf Proposition #1} }
\def \Th #1 {\vs{2mm}\nin {\bf Theorem #1} }
\def \Cor #1 {\vs{2mm}\nin {\bf Corollary #1} }

\def \Def #1 {\vs{2mm}\nin {\bf Definition #1} }
\def \part #1 {\hfil\break\hglue 12pt {\rm (#1)~}}

\def\fs{\footnotesize}

\def\Epi{^{e_\pi}}
\def\Epi1m{^{e_{\pi_1}\mathfrak m}}
\def\epi2m{^{e_{\pi_2}\mathfrak m}}

\begin{document}
\title{
\bf\LARGE  Remarks on two theorems in linear algebra \thanks{This research was supported by the NSF of Shanghai (Grant No. 19ZR1424100). } }
\author{{Yifan Ren\thanks{ivan\_ren@sjtu.edu.cn}}\\
{\small Zhiyuan College,  Shanghai Jiao Tong University}\\
 {Tongsuo Wu\thanks{Corresponding author. tswu@sjtu.edu.cn} }\\
 {\small School of Mathematical Sciences,  Shanghai Jiao Tong University
}
}

\date{}
\maketitle

\begin{center}
\begin{minipage}{12cm}

\vs{3mm}\nin{\small\bf Abstract.} {\fs In this note, we use the concept of a polynomial ring to give a direct proof to Cayley-Hamilton Theorem. We also give an elementary proof to Birkhoff theorem on Bi-stochastic matrices.}

\vs{3mm}\nin {\small Key Words:} {\small polynomial ring; Cayley-Hamilton Theorem; Bi-stochastic matrices; Birkhoff theorem.}

\vs{3mm}\nin {\small 2010 AMS Classification:} {\small Primary: 13P05; 13B25; Secondary: 12D05; 05C50.}

\end{minipage}
\end{center}


Recent rediscovery \cite{2019} on eigenvalues and eigenvectors aroused interests of many readers, and leaves the hope that many more new discoveries can be found even in the traditional field of studies. In this short note, we record two things happened in the classroom discussions.

\section{1. A new proof of Cayley-Hamilton Theorem}

The main result concerning the eigen polynomial $f_1(x)=:{\rm det}(xE-A)$ of a square matrix $A$ is the  Cayley-Hamilton theorem:

\begin{mthm} Let $A$ be a square matrix over a commutative ring $S$. Then the eigen polynomial $f_1(x)=:{\rm det}(xE-A)$ of $A$ annihilates $A$, i.e., $f_1(A)=0$.
\end{mthm}

A typical approach to this theorem is found in many text books, e.g.,  \cite[Page 297]{Shi}, \cite[Theorem 2.5.1]{2001Serre}, or \cite[Theorem 6.4.1]{Guo}. See also \cite[Theorem 2.4.3.2, Ex. 2.4 P3]{HB 2013} for other treatments without using Jordan theorem.

In order to give an alternative and easy proof, recall the definition of a  polynomial ring $R[x]$, where $R$ is any ring with identity. Recall that $R[x]$ is a free left $R$-module with a typical basis $1,x,x^2,\ldots$. Recall that the multiplication in $R[x]$ is defined by $rx=xr$ ($\forall r\in R$) and the distributive law.

\begin{mylem} Let $f(x),g(x),h(x)$ and u(x) be polynomials of $R[x]$ with $f(x)=g(x)h(x)+u(x)$. For any element $r$ of $R$, if $r$ commutes with all coefficients of terms of either $h(x)$ or $g(x)$, then $f(r)=g(r)h(r)+u(r)$ holds in $R$.
\end{mylem}

\begin{prof} Let $a_ix^i$ be a term of $g(x)$ and, $b_jx^j$ a term of $h(x)$. Then we have $(a_ir^i)(b_jr^j)=(a_ib_j)r^{i+j}$, and the result follows by definition of multiplication in $R[x]$.
\end{prof}

\vs{3mm}\noindent{\it Proof to Theorem 1.1.} Let $S$ be any commutative ring with identity, let  $R=M_n(S)$, let $A\in R$ and let $f(x)={\rm det}(xE-A)$. Then $\{Ef(x), Ex-A,(E-A)^*\}\seq R[x]$,  thus Theorem 1.1 follows from $(Ex-A)^*(Ex-A)=Ef(x)$ after applying Lemma 2.1. \quad\quad$\Box$

\vs{3mm}The following is also an immediate consequence of Lemma 1.2:

\begin{myprn} Let $R$ be any ring, and let $r\in R$. Then for any $f(x)\in R[x]$, there exists a polynomial $q(x)\in R[x]$ such that $f(x)=q(x)(x-r)+f(r)$.
\end{myprn}

Note that $xr=rx$ holds in $R[x]$, thus if we accept $x^nr$ as the standard form of a term in $R[x]$, then there is another version of Proposition 1.3. Thus the generalized B\'{e}zout theorem (\cite[Page 81]{2000 Gantmacher}) also follows from Proposition 1.3 (actually, Lemma 1.2).

\vs{3mm}{\fs {\bf Remark.} We hope that this proof could be read by college students. (Actually, the first author is a first-year undergraduate and, the proof is resulted from class discussions.) For this, $S$ could be  regarded as a number field $\mathbb F$, and note that $R[x]$ is exactly the same with the polynomial ring $\mathbb F[x]$
on the number field $\mathbb F$, except that $\mathbb F$ is replaced by a more general ring $R$, e.g., $R=\mathbb Z,\, \mathbb Z[i],\,\mathbb Z[\sqrt{-2}],\, M_n(\mathbb F)$. Note that the key is Lemma 1.2, thus one needs not to know the concept of modules, and the equality $(M_n(\mathbb F))[x]=M_n(\mathbb F[x])$ is essentially not needed to know. What is needed is the concept of a polynomial ring $R[x]$ over a general ring $R$, in which $R=M_n(\mathbb F)$ is needed in the proof.}

\section{2. Bi-stochastic matrices and Birkhoff theorem}

Recall that a column stochastic matrix is a nonnegative matrix, in which the sum of entries of each column is $1$. A matrix is called bi-stochastic, if both $A$ and $A^T$ are  column stochastic.

Recall the following famous theorem of Birkhoff (see e.g., \cite[Theorem 8.7.2]{HB 2013} or \cite[Theorem 5.5.1]{2001Serre}):

\begin{mthm} A nonnegative real matrix $M\in M_n(\mathbb R)$ is bi-stochastic if and only if it is the mass center of some permutation matrices.
\end{mthm}

In order to prove the theorem, techniques from convex geometry is often applied, and one can refer to \cite[Theorem 5.5.1]{2001Serre}) or the first edition of  \cite{HB 2013}. In the second edition of  \cite{HB 2013}, the authors have a noble try to give a more easy and elementary approach to the theorem. While the key is to prove Lemma 8.7.1, by taking advantage of ${\rm det}(xE-A)$. In the proof of the lemma, there is an obstacle that we can not overcome. In the following, we provide an alternative way:

\begin{mylem} {\rm( \cite[Lemma 8.7.1]{HB 2013})} Let $A\in M_n(\mathbb R)$ be a bi-stochastic
matrix, which is not a permutation matrix. Then there is a rearrangement $i_1,\ldots,i_n$ of $1,2,\ldots,n$ such that $a_{k\,i_k}\not=0$ holds for all $k$.
\end{mylem}

\begin{prof} If an entry of $A$ is $1$, then the result follows by induction assumption. In the
following, assume $a_{ij}\not=1,\forall i,j$.

Let \,$a_{i_1\,j_1}\not=0$. Then there exists an integer $j_2\not=j_1$ such that $a_{i_1\,j_2}\not=0$ holds. Clearly, we also have $i_2\not=i_1$ such that $a_{i_2\,j_2}\not=0$.  Again by assumption, there is $j_3\not=j_2$ such that $a_{i_2\,j_3}\not=0$ holds. If further \,$j_3\not=j_1$ holds true, then we have $i_3\not=i_2$ such that
$a_{i_3\,j_3}\not=0$ holds. If further $i_3\not=i_1$, then we continue the process. If $a_{i_3\,j_1}\not=0$, then we finish to obtain
$$\begin{pmatrix}
a_{i_1\,j_1}&a_{i_1\,j_2}&&\\
&a_{i_2\,j_2}&a_{i_2\,j_3}&\\
a_{i_3\,j_1}&&a_{i_3\,j_3}&\\
\end{pmatrix}.$$
On the other hand, if $a_{i_3\,j_1}=0$, then we continue the previous discussion. Clearly, row indexes (and column index) will repeat after finite steps. Without loss of generality, we assume row indexes will repeat first, and assume further that the first repeat appears in $j_{r+1}=j_1$.
This implies that both set $i_1,\ldots,i_r $ and set $j_1,\ldots, j_r$ has cardinality $r$, and
$$a_{i_t\,j_t}\not=0, a_{i_t\,j_{t+1}}\not=0,\,\,\forall t=1,\ldots r.$$
Now we are ready to construct a matrix $B$ in the following:
$$b_{i_t\,j_t}=1, b_{i_t\,j_{t+1}}=-1,\,\,\forall t=1,\ldots r.$$
Now for $c=\underset{1\le t\le r}{{\rm min}}\,\,\{a_{i_t\,j_t},\,\,a_{i_t\,j_{t+1}}\}$,
the matrices $A+cB$  and $A-cB$ are both bi-stochastic, and at least one of them has  zero entries fewer than $A$. Note also that neither has newly added zero entry. This completes the proof by induction.\end{prof}


\begin{thebibliography}{gg}

\bibitem{2019} P.B. Denton, S.J. Parke, T. Tao and X.N Zhang. Eigenvectors from eigenvalues. arXiv:1908.03795v1 [math.RA] 10 Aug 2019.

\bibitem{2000 Gantmacher} F.R. Gantmacher. The Theory of Matrices vol.1. AMS Chelsea Publishing, American Mathematical Society Providence, Rhode Island 2000.

\bibitem{Guo} Y.Q. Guo, K.P. Shum and Z.P. Wang. A Course on Advanced Algebra. Science Press, 2014 Beijing.  (In Chinese)

\bibitem{HB 2013} R.A. Horn and C.R. Johnson.\,\, Matrix Analysis. Cambridge Unversity Press. New York. 2'ed Ed. $2013.$ (1'st Ed. 1985)

\bibitem{2001Serre} D. Serre. \,Matrices: Theory and Applications.  GTM $216$, Springer. New York/Berlin/Heidelberg $2002.$

\bibitem{Shi} S.M. Shi and E.F. Wang. Advanced Algebra.   Higher Education Press, Beijing, Third Ed. 1988. ( In Chinese)


\end{thebibliography}
\end{document}